\documentclass{llncs}
\usepackage[utf8]{inputenc}
\usepackage[T1]{fontenc}
\usepackage{booktabs}
\usepackage{lscape}
\usepackage{cite}
\usepackage{mathtools,amssymb,amsmath}
\usepackage{enumerate}
\usepackage{varioref}
\usepackage{amsmath,mathpazo}
\usepackage[pdfpagelabels=true]{hyperref}
\usepackage{cleveref}
\usepackage{xcolor}
\usepackage{algpseudocode}
\usepackage{algorithm}
\usepackage{listings}
\usepackage{csquotes}
\usepackage{url}
\newif\iflongversion

\hypersetup{
	linktoc = page,
	pdfpagemode = UseNone,
	colorlinks,
	linkcolor={red!50!black},
	citecolor={blue!50!black},
	urlcolor={blue!80!black}
}

\usepackage{paralist}
\usepackage{ragged2e}
\usepackage{tikz}
\usepackage{pgfplots}
\usepackage{tikz-qtree}
\usepackage[linguistics]{forest}
\usetikzlibrary{arrows}
\usetikzlibrary{backgrounds,decorations.pathreplacing, bending}
\usepackage{breqn}
\newcommand*{\F}{\mathbf{F}}

\definecolor{identifiercolor}{rgb}{.4,.6,.56}
\definecolor{stringcolor}{gray}{0.5}
\definecolor{inactivecolor}{rgb}{0.15,0.15,0.5}

\usepackage{listings}
\usepackage[framemethod=TikZ]{mdframed}

\newcounter{theo}[section] 
\renewcommand{\thetheo}{\arabic{section}.\arabic{theo}}
\newenvironment{theo}[2][]{%
\refstepcounter{theo}%
\ifstrempty{#1}%
{\mdfsetup{%
frametitle={%
\tikz[baseline=(current bounding box.east),outer sep=0pt]
\node[anchor=east,rectangle,fill=blue!20]
{\strut Theorem~\thetheo};}}
}%
{\mdfsetup{%
frametitle={%
\tikz[baseline=(current bounding box.east),outer sep=0pt]
\node[anchor=east,rectangle,fill=blue!20]
{\strut Example~\thetheo:~#1};}}%
}%
\mdfsetup{innertopmargin=10pt,linecolor=blue!20,%
linewidth=2pt,topline=true,%
frametitleaboveskip=\dimexpr-\ht\strutbox\relax
}
\begin{mdframed}[]\relax%
\label{#2}}{\end{mdframed}}

\begin{document}

\title{A Note on the Ramanujan Machine}

\author{\'Eric Brier \and David Naccache\inst{1} \and Ofer Yifrach-Stav\inst{1}} 

    \institute{
    DI\'ENS, \'ENS, CNRS, PSL University, Paris, France\\
45 rue d'Ulm, 75230, Paris \textsc{cedex} 05, France\\
\email{\url{ofer.friedman@ens.fr}}, 
\email{\url{david.naccache@ens.fr}}
}

\maketitle      

\begin{abstract}

The Ramanujan Machine project detects new expressions related to constants of interest, such as $\zeta$ function values, $\gamma$ and algebraic numbers (to name a few).

In particular the project lists a number of conjectures involving even and odd $\zeta$ function values, logarithms etc.

We show that many relations detected by the Ramanujan Machine Project stem from a specific algebraic observation and show how to generate infinitely many.

This provides an automated proof and/or an explanation of many of the relations listed as conjectures by the project (although not all of them).
\end{abstract}

\section{Introduction}

The Ramanujan Machine project \cite{rama,ref1,ref2} detects new expressions related to constants of interests, such as $\zeta$ function values, $\gamma$ and various algebraic numbers (to name a few). 

In particular the project lists several of conjectures\footnote{\url{http://www.ramanujanmachine.com/wp-content/uploads/2022/07/results_different_zeta_orders.pdf}} concerning values of the $\zeta$ function.

We show that many of the relations detected by the Ramanujan Machine Project stem from a specific algebraic observation and show how to generate and machine-prove infinitely many.

This provides an automated proof and/or an explanation of many of the relations listed as conjectures by the project (although not all of them).

\section{Theoretical Preamble}

Consider continued fractions defined by the formulae

\begin{align*}
&\left\{
  \begin{aligned}
p_n &= b_n p_{n-1} + a_n p_{n-2}\\
q_n &= b_n q_{n-1} + a_n q_{n-2}\\
  \end{aligned}
\right.
\mbox{~~with the initial conditions:}
&\left\{
  \begin{aligned}
  p_{-1}&=1\\
  p_0&=b_0\\
  \end{aligned}
\right.
\mbox{~~and~~}
&\left\{
  \begin{aligned}
  q_{-1}&=0\\
  q_0&=1\\
  \end{aligned}
\right.
\end{align*}

Let $f$ and $g$ be two functions from which we build $\forall n \geq 1$:

\begin{align*}
\left\{
  \begin{aligned}
 a_n  &= - f(n)^2\\
 b_n  & = \frac{f(n+1) g(n+1) + f(n)  g(n-1)}{g(n)} \\
  \end{aligned}
\right.
\mbox{~~with the initial condition:~~} b_0 = f(1)  g(1)
\end{align*}

We further require $f$ and $g$ to be nonzero for positive integers.

Define the sequence $\F_n$ by:

$$ F_n = \prod_{i=1}^{n+1} f(i)\mbox{~~with the initial conditions:~~} \F_{-1}=1\mbox{~~and~~}\F_0=f(1)$$

We will now prove by induction that 

$$ \forall n \geq -1, p_n = \F_n g(n+1)$$

The above initial conditions ensure this equality for $n=-1$ and $n=0$. Assume that the hypothesis holds true for all values below $n$ and let us compute:

\begin{align*}
  \begin{aligned}
 p_n  &= b_n p_{n-1} + a_n p_{n-2} \\
      &= b_n \F_{n-1} g(n) - f(n)^2 \F_{n-2} g(n-1) \\
      &= \F_{n-1} (f(n+1) g(n+1) + f(n)  g(n-1)) - f(n) \F_{n-1} g(n-1)\\
      &= \F_{n-1} f(n+1) g(n+1)\\
      &= \F_{n} g(n+1)\\
  \end{aligned}
\end{align*}

as desired, thereby proving a closed form for the convergents $p_n$.

We will now handle $q_n$ by first remarking that:

$$ \frac{q_n}{p_n} - \frac{q_{n-1}}{p_{n-1}} = \frac{q_n p_{n-1}-q_{n-1} p_{n}}{p_n p_{n-1}}$$

The recurrence conditions for $p_n$ and $q_n$ ensure that:

$$ q_n p_{n-1}-q_{n-1} p_{n} = -a_n (q_{n-1} p_{n-2}-q_{n-2} p_{n-1}) $$

Induction yields:

$$ q_n p_{n-1}-q_{n-1} p_{n} = (\prod_{i=1}^{n}-a_n) (q_{-1} p_{-2}-q_{-2} p_{-1}) $$

The last term is equal to 1 given the initial conditions. The $n$-term product is, by the definition of the sequence $a_n$, equal to $(\F_{n-1})^2$. 

We hence get:

\begin{align*}
  \begin{aligned}
\frac{q_n}{p_n}-\frac{q_{n-1}}{p_{n-1}}  &=  \frac{(\F_{n-1})^2}{p_n p_{n-1}}\\
      &= \frac{(\F_{n-1})^2}{\F_n g(n+1) \F_{n-1} g(n)}\\
      &= \frac{1}{f(n+1) g(n) g(n+1)}\\
  \end{aligned}
\end{align*}

Once again, by induction (taking into account the initial values), we get:

$$ \frac{q_n}{p_n}  = \sum_{i=0}^{n} \frac{1}{f(i+1) g(i) g(i+1)}$$

The limit $L$ of the continued fraction is thus given by the equation:

$$ \dfrac{1}{L}=\sum_{n=0}^{\infty} \dfrac{1}{f(n+1) g(n) g(n+1)} $$

It suffices now to resort to standard partial fraction decomposition to get relations such as those given by the Ramanujan Machine Project. This is done automatically by symbolic computation software such as Mathematica to avoid tedious yet standard formula manipulation by hand.

\begin{theo}[The Ramanujan Machine identity $\zeta(4)+4\zeta(3)-8$]{thm:pythagoras}
Consider the Ramanujan Project identity\footnote{Note that our notations of $a_n$ and $b_n$ are reversed with respect to theirs} where $a_n=-n^8$ and $b_n=n^4+(n+1)^4+2(n^2+(n+1)^2)$. Posing $g(n)=\eta_1 n+\eta_0$ and identifying we get:

$$
b_n g(n) - ((n + 1)^4  g(n+1) + n^4 g(n - 1)) = (1 + 2 n + 2 n^2) (2 \eta_0 - \eta_1) =0
$$

Which gives the solution $\{\eta_1,\eta_0\}=\{2,1\}$ (the first example processed by our code).
\end{theo}

The last step connecting the observation and the Ramanujan Machine Project is somewhat technical we hence break it down into successive steps: 

\paragraph{Step 1:} The Ramanujan Machine Project considers that the $a_n,b_n$ defining the continued fraction are polynomials. Because nothing in our analysis imposes that the $a_n,b_n$ are polynomials we can write:

{\scriptsize
\begin{equation*}\label{eq:gcf}
    \dfrac{1}{L}=\frac{f(1) g(1)+f(0) g(-1)}{g(0)}-{\cfrac {f(1)^2}{\frac{f(2) g(2)+f(1) g(0)}{g(1)}-{\cfrac {f(2)^2}{\frac{f(3) g(3)+f(2) g(1)}{g(2)}-{\cfrac {f(3)^2}{\frac{f(4) g(4)+f(3) g(2)}{g(3)}-{\cfrac {f(4)^2}{\frac{f(5) g(5)+f(4) g(3)}{g(4)}-\cdots }}}}}}}}
\end{equation*}}

{\scriptsize
\begin{equation*}\label{eq:gcf}
    \dfrac{1}{L}=\frac{f(1) g(1)+f(0) g(-1)}{g(0)}-{\cfrac {g(1)f(1)^2}{\frac{f(2) g(2)+f(1) g(0)}{1}-{\cfrac {g(2)g(1)f(2)^2}{\frac{f(3) g(3)+f(2) g(1)}{1}-{\cfrac {g(3)g(2) f(3)^2}{\frac{f(4) g(4)+f(3) g(2)}{1}-{\cfrac {g(4)g(3)f(4)^2}{\frac{f(5) g(5)+f(4) g(3)}{1}-\cdots }}}}}}}}
\end{equation*}}

{\scriptsize
\begin{equation*}\label{eq:gcf}
    \dfrac{g(0)}{L}=\frac{f(1) g(1)+f(0) g(-1)}{1}-{\cfrac {g(1)g(0)f(1)^2}{\frac{f(2) g(2)+f(1) g(0)}{1}-{\cfrac {g(2)g(1)f(2)^2}{\frac{f(3) g(3)+f(2) g(1)}{1}-{\cfrac {g(3)g(2) f(3)^2}{\frac{f(4) g(4)+f(3) g(2)}{1}-{\cfrac {g(4)g(3)f(4)^2}{\frac{f(5) g(5)+f(4) g(3)}{1}-\cdots }}}}}}}}
\end{equation*}}

{\scriptsize
\begin{equation*}\label{eq:gcf}
   \dfrac{g(0)}{L}={f(1) g(1)+f(0) g(-1)}-{\cfrac {g(1)g(0)f(1)^2}{{f(2) g(2)+f(1) g(0)}-{\cfrac {g(2)g(1)f(2)^2}{{f(3) g(3)+f(2) g(1)}-{\cfrac {g(3)g(2) f(3)^2}{{f(4) g(4)+f(3) g(2)}-\cdots}}}}}}
\end{equation*}}

\paragraph{Step 2:} In the Ramanujan Project many formulae have an $a_n$ (in our notations, i.e. $b_n$ in theirs) of the form $Q=(n+\alpha)^2(n+\beta)(n+\gamma)$, i.e. with a total degree of 4 and a term having a power of two. In many other cases the $a_n$ can be written as $\phi_1(n)^2\phi_2(n)\phi_3(n)$ for some functions $\phi_1,\phi_2,\phi_3$. The form ($f(n)^2g(n-1)g(n)$) at the numerator of the continued fraction obtained in Step 1 explains why:

\paragraph{Step 3:} Operate the change of variable $n_0=n+\alpha$ to get: $$Q={n_0}^2(n_0-\alpha+\beta)(n_0-\alpha+\gamma)$$

\paragraph{Step 4:} Operate the change of variables $\beta_1=-\alpha+\beta$, $\gamma_1=-\alpha+\gamma$ to get: $$Q={n_0}^2(n_0+\beta_1)(n_0+\gamma_1)$$

\paragraph{Step 5:} Operate the change of variable $n_0=n_1(\gamma_1-\beta_1)$ to get:

$$Q={(n_1(\gamma_1-\beta_1))}^2(n_1(\gamma_1-\beta_1)+\beta_1)(n_1(\gamma_1-\beta_1)+\gamma_1)$$

\paragraph{Step 6:} Kick the parasite factor $(\gamma_1-\beta_1)^2$ out of the continued fraction and integrate it in the limit. We get:

$$Q'=\frac{Q}{(\gamma_1-\beta_1)^2}={n_1}^2(n_1(\gamma_1-\beta_1)+\beta_1)(n_1(\gamma_1-\beta_1)+\gamma_1)$$

\paragraph{Step 5:} Declare $f(n_1)=\mbox{ID}(n_1)=n_1$ and $g(n_1)=n_1(\gamma_1-\beta_1)+\beta_1$ to get:

$$Q'=f(n_1)^2g(n_1)(n_1(\gamma_1-\beta_1)+\gamma_1)$$

\paragraph{Step 7:} Observe that:

$$g(n_1-1)=(n_1-1)(\gamma_1-\beta_1)+\gamma_1=(\gamma_1-\beta_1)n_1+\beta_1
$$

and hence:

$$Q'=f(n_1)^2g(n_1)g(n_1-1)$$

\paragraph{Step 8:} Replace in the simplified continued fraction the terms of the form:

$$f(n)g(n)+f(n-1)g(n-2)$$ 

by their variable-changed expression in $n_1$. 

\paragraph{Step 9:} Light an altar candle hoping that the variable changes did not alter the initial conditions. If so a new relation of the form $a_n=(n+\alpha)^2(n+\beta)(n+\gamma)$ and $b_n=\mbox{polynomial}$ was found.

\remark{\textbf{Non Polynomial Rational Fractions: }} Nothing in the preamble assumed that $f$ and $g$ are polynomials or rational fractions, hence any functions satisfying the few properties announced can be used to derive ``magic'' continued fractions provided that there is a way to write the infinite sum under a fancy closed form. We give a few examples in the next section.

\section{Implementation}

The implementation assumes that $f(0)=0$ to enforce the initial conditions\footnote{It is also possible to tweak the code to work with other initial conditions provided that $q_{-1} p_{-2}-q_{-2} p_{-1}=1$, we did not do this here.}, sets:

$$a_n=-f(n)^2\mbox{~~and~~}b_n=\frac{f(n+1) g(n+1)+f(n) g(n-1)}{g(n)}$$

and prints:

$$L=g(0)^2\sum_{i=0}^{\infty}\frac{1}{f(i + 1) g(i)g(i + 1)}$$

As well as the approximate numerical values of both the exact expression and the continued fraction (to visually compare both).

For the sake of compactness we display $L$ and not its inverse.

The code takes $f$ and $g$ from an example list \texttt{Ex} into which the reader can plug any desirable function to generate new relations at wish. In the listing above we changed the order of the printed formulae to fit the longest examples in a landscape layout. In the formulae $C$ stands for Catalan's constant.

\begin{lstlisting}[extendedchars=true,language=Mathematica]

Ex={{z^4,1+2z},{z^3,z+1},{z^7,z+1},{z^4/(z+2),(z+3)},{z^5(z+1),z+1},{z^6/(z+2),(z+1)/(z^2+1)},{z^4(z+1)/(z+2),(z+1)/(z^2+1)},{z^2,(1+Sqrt[2](z^2+z))},{z^4,(1-Sqrt[2](z^2+z))},{z^6,(1+29(z^2+z))},{z^5,(1+4(z^2+z))},{E^(-8-2z)z^3,E^z},{E^(-8-2z)z(2+5z+2z^2)^2,E^z},{E^(-2-2z)z(1+z)^2,E^z},{E^(-2z)z(2+z)^2,E^z},{E^(-2-2z)z(1+4z)^2,E^z},{E^(-2z)z(z+z^2)^2,E^z},{E^(-2z)z(1+2z+z^2)^2,E^z},{E^(-2-2z)z(z+2z^2)^2,E^z},{E^(-2z)z(3z+4z^2)^2,E^z},{9z/Exp[z],Exp[z]},{(-2-3z)z,(2+z)^10},{(-2-2z)z,(2+z)^10},{(-3-z)z,(2+z)^10},{-2z,(2+z)^10},{z(-2+3z),(2+z)^10},{3z^2,(2+z)^10},{z(1+3z),(2+z)^10},{(2z+1)/E^z,E^z},{(2z+3)/E^z,E^z},{(2z^2+z)/E^z,E^z},{(2z^2+3z+1)/E^z,E^z},{z^2,1+z},{z(1+2z),1+z},{z(1+3z),1+z},{z(2+z^2),2+z},{z(2+z+z^2),3+z},{z(2+2z+z^2),4+z},{z(2+4z+z^2),3+z},{z(4+2z^2),2+z},{z(z+3z^2),1+z},{z(2z+3z^2),1+z},{z(4z+3z^2),1+z},{z(3z+4z^2),1+z},{z(3z+4z^2),2+z}};

F:=Function[{w,x},x[[1]]/.{z->w}]; 
G:=Function[{w,x},x[[2]]/.{z->w}]; 
For[i=1,i<=Length[Ex],
 v=.; 
 CF=(F[v+1,Ex[[i]]] G[v+1,Ex[[i]]]+F[v,Ex[[i]]] G[v-1,Ex[[i]]])/G[v,Ex[[i]]];
 st=N[{1,0,CF/.v->0,1},10000];
 For[n=1,n<= 1000,n++,
  bn=CF/.v->n;
  an=-F[n,Ex[[i]]]^2;
  st={{0,0,1,0},{0,0,0,1},{an,0,bn,0},{0,an,0,bn}}.st];
 formalexpr=G[0,Ex[[i]]]^2/F[t+1,Ex[[i]]]/G[t,Ex[[i]]]/G[t+1,Ex[[i]]];
 closedform=Simplify[Sum[formalexpr, {t,0,Infinity}]];
 approxform=N[st[[4]]/st[[3]],1000];
 Print[{closedform,N[closedform,100],N[approxform,100]}];
 i++]

\end{lstlisting}

{\tiny
$$\left\{v^4,2 v+1,8-\frac{2 \pi ^2}{3}-\frac{\pi ^4}{90}\right\}$$

$$\left\{v^3,v+1,-\zeta (3)+\frac{\pi ^4}{90}+\frac{\pi ^2}{6}-1\right\}$$

$$\left\{v^7,v+1,-\zeta (3)-\zeta (5)-\zeta (7)+\frac{\pi ^8}{9450}+\frac{\pi ^6}{945}+\frac{\pi ^4}{90}+\frac{\pi ^2}{6}-1\right\}$$

$$\left\{\frac{v^4}{v+2},v+3,\frac{1}{270} \left(-270 \zeta (3)+9 \pi ^4+15 \pi ^2-55\right)\right\}$$

$$\left\{v^5 (v+1),v+1,-4 \zeta (3)-2 \zeta (5)+\frac{\pi ^6}{945}+\frac{\pi ^4}{30}+\pi ^2-7\right\}$$

$$\left\{\frac{v^6}{v+2},\frac{v+1}{v^2+1},11 \zeta (3)+10 \zeta (5)+4 \zeta (7)-\frac{5 \pi ^2}{3}-\frac{11 \pi ^4}{90}-\frac{2 \pi ^6}{315}+10\right\}$$

$$\left\{\frac{v^4 (v+1)}{v+2},\frac{v+1}{v^2+1},4 (5 \zeta (3)+\zeta (5)+13)-\frac{41 \pi ^2}{6}-\frac{\pi ^4}{9}\right\}$$

$$\left\{v^2,\sqrt{2} \left(v^2+v\right)+1,\frac{\pi ^2 \left(\sqrt{2}-4\right)-30 \sqrt{2}+36+3 \pi  \left(3 \sqrt{2}-2\right) \sqrt{2 \sqrt{2}-1} \tanh \left(\frac{1}{2} \sqrt{2 \sqrt{2}-1} \pi \right)}{6 \left(\sqrt{2}-4\right)}\right\}$$

$$\left\{v^4,1-\sqrt{2} \left(v^2+v\right),-4 \sqrt{2}+\frac{\pi ^4}{90}+\frac{\pi ^2 \left(3 \sqrt{2}+5\right)}{6 \sqrt{2}+3}-5-\frac{\left(27 \sqrt{2}+38\right) \pi  \tan \left(\frac{1}{2} \sqrt{2 \sqrt{2}+1} \pi \right)}{2 \left(2 \sqrt{2}+1\right)^{3/2}}\right\}$$

$$\left\{v \left(v^2+2 v+2\right),v+4,\frac{49}{18}-\frac{4}{5} \pi  \coth (\pi )\right\}$$

$$\left\{v \left(v^2+4 v+2\right),v+3,\frac{1}{8} \left(9 \sqrt{2} \pi  \cot \left(\sqrt{2} \pi \right)-10\right)\right\}$$

$$\left\{v \left(2 v^2+4\right),v+2,\frac{1}{12} \left(6-\sqrt{2} \pi  \coth \left(\sqrt{2} \pi \right)\right)\right\}$$

$$\left\{v \left(3 v^2+v\right),v+1,\zeta (3)-\frac{2 \pi ^2}{3}-\frac{9 \sqrt{3} \pi }{4}+40+\log \left(\frac{1}{3486784401 \sqrt[4]{3}}\right)\right\}$$

$$\left\{v \left(3 v^2+2 v\right),v+1,\frac{1}{48} \left(-3 (-8 \zeta (3)-65+\log (3)+16 \log (243))+27 \pi  \sqrt{3}-10 \pi ^2\right)\right\}$$

$$\left\{v \left(3 v^2+4 v\right),v+1,\frac{1}{768} \left(192 \zeta (3)-56 \pi ^2+54 \pi  \sqrt{3}-447-324 \log (2)+162 \log (3)+324 \log (6)\right)\right\}$$

$$\left\{v \left(4 v^2+3 v\right),v+1,\frac{1}{162} \left(54 \zeta (3)-21 \pi ^2+192 \pi +350-384 \log (8)\right)\right\}$$

$$\left\{v \left(4 v^2+3 v\right),v+2,\frac{1}{270} \left(30 \pi ^2-768 \pi -1049+1536 \log (8)\right)\right\}$$

$$\left\{e^{-v} \left(2 v^2+3 v+1\right),e^v,e \left(-1+e (\log (e-1)-1)+2 \sqrt{e} \tanh ^{-1}\left(\frac{1}{\sqrt{e}}\right)\right)\right\}$$

$$\left\{v^2,v+1,\zeta (3)-\frac{\pi ^2}{6}+1\right\}$$

$$\left\{v (2 v+1),v+1,\frac{\pi ^2}{6}-7+\log (256)\right\}$$

$$\left\{v (3 v+1),v+1,\frac{1}{12} \left(2 \pi ^2+9 \pi  \sqrt{3}-156+81 \log (3)\right)\right\}$$

$$\left\{v \left(v^2+2\right),v+2,1-\frac{\pi  \coth \left(\sqrt{2} \pi \right)}{3 \sqrt{2}}\right\}$$

$$\left\{v \left(v^2+v+2\right),v+3,\frac{25}{16}-\frac{9 \pi  \tanh \left(\frac{\sqrt{7} \pi }{2}\right)}{8 \sqrt{7}}\right\}$$

$$\left\{e^{-v} (2 v+1),e^v,e \left(\sqrt{e} \tanh ^{-1}\left(\frac{1}{\sqrt{e}}\right)-1\right)\right\}$$

$$\left\{e^{-v} (2 v+3),e^v,\frac{1}{3} e \left(-3 e-1+3 e^{3/2} \tanh ^{-1}\left(\frac{1}{\sqrt{e}}\right)\right)\right\}$$

$$\left\{e^{-v} \left(2 v^2+v\right),e^v,-e \left(-3+\log (e-1)+2 \sqrt{e} \tanh ^{-1}\left(\frac{1}{\sqrt{e}}\right)\right)\right\}$$

$$\left\{v^6,29 \left(v^2+v\right)+1,\frac{\pi ^6}{945}+\frac{87 \pi ^4}{10}+\frac{306124 \pi ^2}{3}-\frac{478731681}{2}+\frac{15895741}{10} \pi  \sqrt{29} \tan \left(\frac{5 \pi }{2 \sqrt{29}}\right)\right\}$$

$$\left\{v^5,4 \left(v^2+v\right)+1,8 \zeta (3)+\zeta (5)+56-48 \log (4)\right\}$$

$$\left\{e^{-2 v-8} v^3,e^v,e^9 \zeta (3)\right\}$$

$$\left\{e^{-2 v-8} v \left(2 v^2+5 v+2\right)^2,e^v,-\frac{1}{108} e^9 \left(13 \pi ^2+8 (\log (16)-19)\right)\right\}$$

$$\left\{e^{-2 v-2} v (v+1)^2,e^v,-\frac{1}{6} e^3 \left(\pi ^2-12\right)\right\}$$

$$\left\{e^{-2 v} v (v+2)^2,e^v,e-\frac{e \pi ^2}{12}\right\}$$

$$\left\{e^{-2 v-2} v (4 v+1)^2,e^v,-\frac{1}{4} e^3 \left(8 C+\pi ^2+2 \pi +4 (\log (8)-8)\right)\right\}$$

$$\left\{e^{-2 v} v \left(v^2+v\right)^2,e^v,e \left(\zeta (3)-\frac{\pi ^2}{2}+4\right)\right\}$$

$$\left\{e^{-2 v} v \left(v^2+2 v+1\right)^2,e^v,-\frac{1}{90} e \left(90 (\zeta (3)-4)+\pi ^4+15 \pi ^2\right)\right\}$$

$$\left\{e^{-2 v-2} v \left(2 v^2+v\right)^2,e^v,-\frac{1}{3} e^3 \left(-3 (\zeta (3)+32)+5 \pi ^2+36 \log (4)\right)\right\}$$

$$\left\{e^{-2 v} v \left(4 v^2+3 v\right)^2,e^v,-\frac{1}{243} e \left(-288 C-27 \zeta (3)+48 \pi ^2-72 \pi -256+144 \log (8)\right)\right\}$$

$$\left\{9 e^{-v} v,e^v,\frac{1}{9} (e-e \log (e-1))\right\}$$}

\begin{landscape}
{\tiny
$$\left\{v^4,\sqrt{2} \left(v^2+v\right)+2,\frac{4 \pi ^4 \left(\sqrt{2}-8\right)+30 \pi ^2 \left(17 \sqrt{2}-12\right)+45 \left(64-39 \sqrt{2}\right)+90 \pi  \left(5 \sqrt{2}-11\right) \sqrt{4 \sqrt{2}-1} \tanh \left(\frac{1}{2} \sqrt{4 \sqrt{2}-1} \pi \right)}{360 \left(\sqrt{2}-8\right)}\right\}$$

$$\left\{v^6,\sqrt{2} \left(v^2+v\right)+2,\frac{16 \pi ^6 \left(\sqrt{2}-8\right)+84 \pi ^4 \left(17 \sqrt{2}-12\right)-6615 \left(19 \sqrt{2}-28\right)+630 \pi ^2 \left(39 \sqrt{2}-64\right)+1890 \pi  \left(21 \sqrt{2}-29\right) \sqrt{4 \sqrt{2}-1} \tanh \left(\frac{1}{2} \sqrt{4 \sqrt{2}-1} \pi \right)}{15120 \left(\sqrt{2}-8\right)}\right\}$$

$$\left\{v^8,\sqrt{2} \left(v^2+v\right)+2,\frac{16 \pi ^8 \left(\sqrt{2}-8\right)+80 \pi ^6 \left(17 \sqrt{2}-12\right)+22050 \pi ^2 \left(19 \sqrt{2}-28\right)+420 \pi ^4 \left(39 \sqrt{2}-64\right)-14175 \left(149 \sqrt{2}-200\right)+9450 \pi  \left(69 \sqrt{2}-91\right) \sqrt{4 \sqrt{2}-1} \tanh \left(\frac{1}{2} \sqrt{4 \sqrt{2}-1} \pi \right)}{151200 \left(\sqrt{2}-8\right)}\right\}$$

$$\left\{v^{10},1-\sqrt{2} \left(v^2+v\right),\frac{2 \pi ^6 \left(8 \sqrt{2}+11\right)}{135 \left(4 \sqrt{2}+9\right)}+\frac{\pi ^{10}}{93555}+\frac{\pi ^8 \left(13 \sqrt{2}+17\right)}{4725 \left(4 \sqrt{2}+9\right)}+\frac{8 \pi ^4 \left(30 \sqrt{2}+43\right)}{45 \left(4 \sqrt{2}+9\right)}-4 \left(79 \sqrt{2}+111\right)+\frac{2 \pi ^2 \left(176 \sqrt{2}+249\right)}{4 \sqrt{2}+9}-\frac{2 \left(2159 \sqrt{2}+3057\right) \pi  \tan \left(\frac{1}{2} \sqrt{2 \sqrt{2}+1} \pi \right)}{\left(2 \sqrt{2}+1\right)^{5/2}}\right\}$$

$$\left\{v^{12},\sqrt{2} \left(v^2+v\right)+2,\frac{88448 \pi ^{12} \left(\sqrt{2}-8\right)+436800 \pi ^{10} \left(17 \sqrt{2}-12\right)+75675600 \pi ^6 \left(19 \sqrt{2}-28\right)+2162160 \pi ^8 \left(39 \sqrt{2}-64\right)+170270100 \pi ^4 \left(149 \sqrt{2}-200\right)+425675250 \pi ^2 \left(1381 \sqrt{2}-1996\right)-638512875 \left(4479 \sqrt{2}-6320\right)+1277025750 \pi  \left(689 \sqrt{2}-963\right) \sqrt{4 \sqrt{2}-1} \tanh \left(\frac{1}{2} \sqrt{4 \sqrt{2}-1} \pi \right)}{81729648000 \left(\sqrt{2}-8\right)}\right\}$$

$$\left\{(-3 v-2) v,(v+2)^{10},\frac{3}{4} (654616480 \zeta (3)+342927872 \zeta (5)+73113600 \zeta (7)+6422528 \zeta (9)+31737050334-2324522934 \log (2)+1162261467 \log (3)+2324522934 \log (6))-\frac{1}{4} 1162261467 \pi  \sqrt{3}-2279750473 \pi ^2-\frac{156065456 \pi ^4}{5}-\frac{913408 \pi ^8}{225}-\frac{141308416 \pi ^6}{315}-\frac{131072 \pi ^{10}}{10395}\right\}$$

$$\left\{(-2 v-2) v,(v+2)^{10},-\frac{512 \left(-467775 (4341764 \zeta (3)+1794064 \zeta (5)+389184 \zeta (7)+41216 \zeta (9)+1024 \zeta (11)+110374897)+48640 \pi ^{10}+10638144 \pi ^8+1035033120 \pi ^6+68103341460 \pi ^4+4803747223275 \pi ^2\right)}{467775}\right\}$$

$$\left\{(-v-3) v,(v+2)^{10},\frac{512 \left(-39106972899225+3428161436700 \pi ^2+47204776080 \pi ^4+649503360 \pi ^6+5144832 \pi ^8+10240 \pi ^{10}\right)}{467775}\right\}$$

$$\left\{v (3 v-2),(v+2)^{10},-\frac{3 (2240882559523414400 \zeta (3)+902645353815040000 \zeta (5)+178178187264000000 \zeta (7)+13985382400000000 \zeta (9)+126999855214834440651+2324522934 \log (2)-1162261467 \log (3)-2324522934 \log (6))}{40000000000}+\frac{1162261467 \pi  \sqrt{3}}{40000000000}+\frac{65536 \pi ^{10}}{22275}+\frac{88928768 \pi ^8}{65625}+c\right\}$$

$$
\mbox{where:~~}c=\frac{162966967472 \pi ^6}{984375}+\frac{111360024325867 \pi ^4}{9375000}+\frac{859444691048496267 \pi ^2}{1000000000}
$$

$$\left\{3 v^2,(v+2)^{10},\frac{256 \left(-467775 (3325728 \zeta (3)+1342656 \zeta (5)+267264 \zeta (7)+21504 \zeta (9)+159629723)+25600 \pi ^{10}+11037312 \pi ^8+1318680000 \pi ^6+93801403080 \pi ^4+6769898635500 \pi ^2\right)}{1403325}\right\}$$

$$\left\{v (3 v+1),(v+2)^{10},-\frac{512 \left(-18711 (-2 (1307738292800 \zeta (3)+529342480000 \zeta (5)+106368000000 \zeta (7)+8800000000 \zeta (9)+55119052506237+1162261467 \log (6))+2324522934 \log (2)-1162261467 \log (3))+7249024769679 \pi  \sqrt{3}+c'\right)}{60908203125}\right\}$$

$$\mbox{where:~~}c'=-800000000 \pi ^{10}-317671200000 \pi ^8-37028146320000 \pi ^6-2609348019339600 \pi ^4-187811501659121340 \pi ^2$$

$$\left\{-2 v,(v+2)^{10},-3072 (106154 \zeta (3)+42664 \zeta (5)+8352 \zeta (7)+640 \zeta (9)+10269357)+\frac{262144 \pi ^{10}}{31185}+\frac{950272 \pi ^8}{225}+\frac{18661376 \pi ^6}{35}+\frac{581261312 \pi ^4}{15}+2814203392 \pi ^2\right\}$$}
\end{landscape}

\section{Using Formal Manipulations}

Mathematica is powerful enough to take into account parametric summations. If we remove from our code the numerical part and add parameters to the input we can get general forms in which we later set parameters at wish.

In the following example we arbitrarily broke the result into several pieces to fit in the page:

\begin{lstlisting}[extendedchars=true,language=Mathematica]
Ex={{((a z)^4+(c z)^3)u,d+b(2+2z)}};
F:=Function[{w,x},x[[1]]/.{z->w}];
G:=Function[{w,x},x[[2]]/.{z->w}];
For[i=1,i<=Length[Ex],v=.;
CF=(F[v+1,Ex[[i]]]G[v+1,Ex[[i]]]+
F[v,Ex[[i]]]G[v-1,Ex[[i]]])/G[v,Ex[[i]]];
formalexpr=
G[0,Ex[[i]]]^2/F[t+1,Ex[[i]]]/G[t,Ex[[i]]]/G[t+1,Ex[[i]]];
closedform=Simplify[Sum[formalexpr,{t,0,Infinity}]];
Print[{closedform}];
i++]
\end{lstlisting}

{\tiny

$$\mbox{result}=\frac{\ell_1+\ell_2+\ell_3+\ell_4+\ell_5+\ell_7+\ell_6}{6 c^9 d^3 u (2 b+d)^2 \left(2 b c^3-a^4 d\right) \left(2 b c^3-a^4 (2 b+d)\right)}$$

$$\ell_1=96 b^3 c^9 (2 b+d) \psi ^{(0)}\left(\frac{d}{2 b}+1\right) \left(b c^3 \left(4 b^2+6 b d+3 d^2\right)-2 a^4 \left(2 b^3+4 b^2 d+3 b d^2+d^3\right)\right)$$

$$
\ell_2=+6 a^{16} d^3 (2 b+d)^4 \psi ^{(0)}\left(\frac{c^3}{a^4}+1\right)+\left(2 b c^3-a^4 d\right) 
$$

$$
\ell_3=384 \gamma  b^6 c^6 \left(c^3-a^4\right)-32 b^5 c^3 d \left(6 \gamma  a^8+\left(24 \gamma -\pi ^2\right) a^4 c^3+\left(\pi ^2-24 \gamma \right) c^6\right)+a^4 d^6 \left(-6 \gamma  a^8+\pi ^2 a^4 c^3-6 c^6 \zeta (3)\right)
$$

$$\ell_4=-16 b^4 d^2 \left(6 \gamma  a^{12}+\left(24 \gamma -\pi ^2\right) a^8 c^3+2 a^4 c^6 \left(3 \zeta (3)-2 \pi ^2+18 \gamma \right)-2 c^9 \left(3 \zeta (3)-2 \pi ^2+18 \gamma \right)\right)$$

$$\ell_5=8 b^3 d^3 \left(-24 \gamma  a^{12}+4 \left(\pi ^2-9 \gamma \right) a^8 c^3-6 a^4 c^6 \left(4 \zeta (3)-\pi ^2+4 \gamma \right)+c^9 \left(18 \zeta (3)-5 \pi ^2+18 \gamma -6\right)\right)$$

$$
\ell_6=-2 b d^5 \left(24 \gamma  a^{12}+2 \left(3 \gamma -2 \pi ^2\right) a^8 c^3-a^4 c^6 \left(\pi ^2-24 \zeta (3)\right)-6 c^9 \zeta (3)\right)
$$

$$
\ell_7=-8 b^2 d^4 \left(18 \gamma  a^{12}+3 \left(4 \gamma -\pi ^2\right) a^8 c^3+a^4 c^6 \left(18 \zeta (3)-2 \pi ^2+3 \gamma \right)+c^9 \left(\pi ^2-9 \zeta (3)\right)\right)
$$}

This is much more practical than the blind exploration of values. Whilst the code does nothing with \texttt{CF}, we left it in the program to explicit the continued fraction being calculated. The chosen example has nothing fundamental or conceptual in it and was chosen at random for illustrative purposes.

\section{Conclusion \& Further Challenges}

Given the above it appears that several of the Ramanujan Project conjectures\footnote{We did not exhaust all the relations listed online.} can be explained and/or automatically machine-proved with virtually no effort.\smallskip 

We did not machine-prove the online conjectures one by one as this implies retyping their polynomials, a retro-solving for $f$ and $g$ by identification and computing the closed forms using formal summation to infinity. Nonetheless, many such relations can be generated automatically quasi-instantaneously on a simple PC. 

An interesting question is that of reversing continued fractions from a target constant. For instance, determine $a_n$ and $b_n$ such that the continued fraction converges to a constant chosen \textsl{a priori}, e.g.:

$$
\frac{1}{L}=\sum_{i=2}^{100}\zeta(i)
$$

We did not research this challenge and leave it to readers interesting in pursuing this line of investigation.

\bibliographystyle{abbrv}
\bibliography{biblio.bib}
\end{document}